\documentstyle[11pt]{article}

\title{Sur une classe de groupoides riemanniens\footnote
{Supported in part by a CNCSU Grant no. 132/95 and 
published in Analele Univ. Timisoara, Ser. Mat.-Inf., Vol. XXXIV, 1, 1996,
137-150
}}

\author{Paul Popescu\\
Departement de mathematiques, Universit\'e de Craiova,\\
13, A.I.Cuza street, Craiova, 1300, Romania}
\date{}
     \newcommand{\btu}{B_T^1}
     \newcommand{\btuu}{$\btu \; $}

     \newcommand{\ftbu}{\overline {F}_T^1}
     \newcommand{\ftbuu}{$\ftbu \; $}
     \newcommand{\frtu}{{\cal F}^1_T}
     \newcommand{\frtuu}{$\frtu \; $}
     \newcommand{\wtu}{W^1_T}
     \newcommand{\wtuu}{$\wtu \; $}
     \newcommand{\fbtu}{\overline {\cal F}^1_T}
     
     \newcommand{\fb}{\overline {\cal F}}
     \newcommand{\mf}{$(M,{\cal F})\; $}
     \newcommand{\bff}{$(\btu,\frtu )\; $}
     \newcommand{\mfb}{$(M,\fb)\; $}
     \newcommand{\bffb}{$(\btu ,\fbtu)\; $}
     \newcommand{\prfib}[1]{\times _{#1}}

     \newcommand{\hfbtu}{{\cal H}(\fbtu )}
     \newcommand{\hfbtuu}{$\hfbtu \; $}
     \newcommand {\grupoid}[4]{{#1}     
\left . \begin{array}  {c}    \vspace {-2  mm}      {#3} 
\\{\longrightarrow} \vspace {-3 mm} \\
 {\longrightarrow}\vspace {-2 mm}\\{#4} \end{array} \right . 
 {#2}}
     \newcommand {\gab}{                
$\grupoid {\Gamma}{\Gamma _0}{\alpha }{\beta }\; $}
     \newcommand {\gabm}{                
$\grupoid {\Gamma}{M}{\alpha }{\beta }\; $}
     \newcommand {\rref}[1]{(\ref {#1})} 
     \newcommand {\ggamma}{$\Gamma \; $} 
     \newcommand {\ggammap}{$\Gamma '\; $} 
     \newcommand {\ggz}{$\Gamma _ 0\; $} 
     \newcommand {\aalpha}{$\alpha \; $} 
     \newcommand {\bbeta}{$\beta \; $}   
        \newcommand{\dem}{\hspace {2 cm}\underline {\em 
                  D\'{e}monstration} } 
        \newcommand{\qed}{{\em (Q.e.d.)}}

     \newtheorem{th}{Th\'{e}or\`{e}me}
     \newtheorem{pr}{Proposition}
     \newtheorem{lm}{Lemme}

\begin{document}
\maketitle 
\begin{abstract}
Dans ce travail on montre qu'il existe un groupoide riemannien dont les orbites sont les adh\'erences des feuilles d'un feuilletage riemannien r\'egulier sur une vari\'et\'e compacte. Ce groupoide est \'equivalent (au sens g\'en\'eralis\'e de Haefliger) avec un groupoide transformationel sur la vari\'et\'e basique.
\end{abstract}
\section{Introduction}
     Nous  allons  utiliser  les  sources   suivantes   pour   les 
d\'{e}finitions,  les  r\'{e}sultats  fondamentaux  et  les  notations  de 
base:  \cite {Mo}  pour  les   Feuiletages   Riemanniens,   \cite 
{AD}  pour  les 
Groupoides de Lie,  \cite {Ha}  pour  l'\'{e}quivalence  des  Groupoides  et 
\cite {GGHR} pour les Groupoides Riemanniens.

     Soit \mf  un  feuilletage  riemannien  (r\'{e}gulier et   de 
codimension $q\; $) sur  une  vari\'{e}t\'{e}  compacte $M\; $ et 
\bff   le feuilletage  \'{e}lev\'{e}  dans  le  fibr\'{e}   
principal   des rep\`{e}res transverses. 
Ce  dernier  a   l'adh\'{e}rence   r\'{e}guli\`{e}re,   donc   le 
feuilletage \bffb est r\'{e}gulier . Il est  aussi  simple,  
d\'{e}fini   par  la  fibration basique $p'\! :\! \btu 
\rightarrow \wtu $  avec la fibre  type \ftbuu   et  la  base  la 
vari\'{e}t\'{e} basique $\wtu \; $ . Il existe une action \`{a} droite
$$\wtu \times O(q)\rightarrow \wtu
$$    
Cette action et  l'action  \`{a}  droite  de  $O(q)\;  $  sur  le 
fibr\'{e}  principal \btuu commutent avec  la  projection  basique  
$p'\! :\! \btu \rightarrow \wtu \; $  et  l'espace  des orbites 
de l'action est hom\'{e}omorphe avec l'espace des feuilles du 
feuilletage singulier \mfb (\cite {Mo}).

      Pour  ce  qui  suit  on  a   besoin   de   certaines   m\'{e}triques 
riemanniennes sur \btuu et \wtuu   .  Pour \c{c}a  on  peut  \'{e}lever  une 
m\'{e}trique adapt\'{e}e au feuilletage riemannien \mf au \btuu   en 
utilisant son parall\'{e}lisme transvers. Cette   m\'{e}trique   est  
\'{e}galement 
adapt\'{e}e aux feuilleteges riemanniens \bff et \bffb .Le 
dernier feuilletage  a  une  holonomie  nulle  et  il  induit  une 
m\'{e}trique riemannienne sur \wtuu  .

     Conform\'{e}ment \`{a} \cite {GGHR}, un groupoide de  Lie    
\gab   est 
riemannien si on a des structures riemanniennes  sur  \ggamma  et 
\ggz ,  tel
que, si \aalpha et \bbeta sont des submersions riemanniennes et  l'inversion 
du groupoide est une isom\'{e}trie. En utilisant cette d\'{e}finition  et 
conform\'{e}ment aux r\'{e}sultats de \cite {Wi} on peut  dire  que  le  
groupoide 
d'holonomie ${\cal H}(\frtu )\; $ est un groupoide riemannien s\'{e}par\'{e}.

     On  peut  d\'{e}finir  une  action  naturelle :
$$ {\cal H}(\frtu )\times O(q)\rightarrow {\cal H}(\frtu )
$$
En factorisant cette action  on  obtient  un  morphisme  des 
groupoides de Lie
\begin{equation}\label {star}
\begin{picture}(165,65)(5,5) 
\put (2,10){${\cal H}(\overline {\cal F}_T^1)$} 
     \put (20,50){$\Gamma $}  \put  (60,10){$B_T^1$} 
      \put (60,50){$M$}  
\put (24,20){\vector (0,1){26}} \put (64,20){\vector (0,1){26}} 
\put (36,10){\vector (1,0){22}} \put (30,50){\vector (1,0){28}} 
\put (36,16){\vector (1,0){22}} \put (30,56){\vector (1,0){28}} 
\put (40,18){$\alpha '$} \put (40,58){$\alpha $}
\put (40,1){$\beta '$} \put (40,41){$\beta $}
\put (26,30){$\pi $} \put (66,30){$p$} 
\end{picture}
\end{equation}

     Le groupoide \gabm  est compact, s\'{e}par\'{e} et \aalpha  , 
\bbeta  sont  des fibrations avec la fibre type \ftbuu .

     On  peut  \'{e}galement  introduire  par  factorisation  une  
m\'{e}trique 
riemannienne sur \ggamma , qui est aussi un groupoide  riemannien.  Dans 
le  diagramme \rref {star}  toutes  les  fl\`{e}ches  sont  des   submersions 
riemanniennes.

     Les orbites \cite {AD} de \ggamma sont les adh\'{e}rences  des  
feuilles  de \mf et on dit que le  feuilletage \mfb  est  
defini   par  le groupoide \ggamma \cite {GGHR} . On a le r\'{e}sultat 
suivant:
\begin{th}\label {thunu}
Soit \mf un feuilletage riemannien r\'{e}gulier sur  la  
vari\'{e}t\'{e}  compacte  $M\; $.   Alors   il   existe   un   
groupoide riemannien s\'{e}par\'{e} et compact \gabm qui 
d\'{e}finisse le feuilletage des  adh\'{e}rences  \mfb .
\end{th}
     En particulier, si on utilise  \cite {GGHR}  on  retrouve  que  le 
feuilletage \mfb est un feuilletage riemannien singulier au sens 
de \cite {Mo}.
     Le {\em Th\'{e}or\`{e}me \ref {thunu}}  donne  une   
r\'{e}ponse   positive \`{a}   quelques probl\`{e}mes    
pos\'{e}s en \cite {GGHR} relativement aux   groupoides 
riemanniens:
\begin{description}
\item[-]Si \mf est un feuilletage transversalement complet  son 
adh\'{e}rence \mfb est r\'{e}gulier  et  un  feuilletage  simple  
\cite {Mo}  donc un   cas   particuli\`{e}re   d'une    fibration    
de    Seifert g\'{e}n\'{e}ralis\'{e}e (\cite {Gh}). 
Mais de \cite {Gh}  un  tel  feuilletage  est  riemannien  et  avec  le 
Th\'{e}or\`{e}me \ref {thunu} on d\'{e}duit qu'il existe un groupoide 
Riemannien qui  le d\'{e}finisse.
\item [-]Le probl\`{e}me de construire un groupoide riemannien  
qui d\'{e}finisse un feuilletage riemannien singulier a une 
solution dans le cas o\'{u} le feuilletage singulier est 
l'adh\'{e}rence  des  feuilles d'un feuilletage riemannien 
r\'{e}gulier.
\end{description}

     Dans le cas g\'{e}n\'{e}ral on ne peut pas esp\'{e}rer que le 
groupoide soit compact ou horisontalement compl\`{e}t  (v.  \cite 
{P}). On  peut  donner  l'interpretation  suivante  du  groupoide 
\ggamma :
\begin{th}\label {thdoi}
Le groupoide transformationel 
$$\Gamma '=\grupoid {\wtu \times O(q)}{\wtu }{ }{ }
$$
est riemannien, les fl\`{e}ches du diagramme
\begin{equation}
\begin{picture}(165,65)(5,5) 
\put (63,10){$M$} \put (60,50){$W_1^T$}  \put  (110,10){$\Gamma $} 
      \put (110,50){$\Gamma '$}  
     \put (12,30){$B_T^1$}
\put (108,10){\vector (-1,0){28}} \put (108,50){\vector (-1,0){28}} 
\put (108,16){\vector (-1,0){28}} \put (108,56){\vector (-1,0){28}} 
\put (27,35){\vector (2,1){30}}
\put (27,28){\vector (2,-1){30}}
\put (90,18){$\alpha $} \put (90,58){$\alpha '$}
\put (90,1){$\beta $} \put (90,41){$\beta '$}
\end{picture}
\end{equation}

sont des submersions riemanniennes et \btuu  est une \'{e}quivalence  des 
groupoides \ggamma et \ggammap ( au sens de \cite {Ha} ).
\end{th} 

     On dit alors que \btuu est  une  \'{e}quivalence  riemannienne  de 
groupoides. 

     Le {\em Th\'{e}or\`{e}me \ref {thdoi}} donne une  
interpr\'{e}tation   de  la m\'{e}thode de P.Molino d'\'{e}tudier 
les feuilletages  riemanniens  (\cite  {Mo}) \`{a} partir  du 
fibr\'{e} des rep\`{e}res  transverses.  On  peut  consid\'{e}rer  aussi  les 
r\'{e}sultats de \cite {Ha} pour \'{e}tudier les implications  
homotopiques  de l'\'{e}quivalence   \'{e}tablie   (par   exemple   
les    classes caract\'{e}ristiques). 
Nous allons \'{e}tudier cela dans un autre travail.

\section{Le groupoide transformationel \ggammap }
     Nous  allons  utiliser  dans  ce  paragraph  \cite {Ha}  pour  les 
fibr\'{e}s G-principaux et \cite {Mo} pour  la  structure  de  feuilletages 
riemanniens.
     L'action \`{a}  droite   de $O(q)\; $   sur \wtuu    d\'{e}finit   le  
groupoide transformationel
$$\Gamma '=\grupoid {\wtu \times O(q)}{\wtu }{\alpha '}{\beta '}
$$
o\'{u} $a'(w,g) = w\cdot g\; ,\; b'(w,g) = w\; ,\;  i(w) = (w,e)\; $ (le  
plongement  de \wtuu dans \ggammap) et $(w,g)  = (w\cdot g , g )$.  On  
a  $(w_1  ,g_1  )\cdot (w_2  ,g_2  )= (w_2 , g_1 \cdot g_2 )\; $ qui est 
d\'{e}finit pour $w_2 =w_1 \cdot  g_1\; $  .
\begin{pr}\label {pru} Relativement aux projections  $  p'\!  :\! 
\btu \rightarrow \wtu $ et $p \! :\! \btu \rightarrow M$ , \btuu  
est un fibr\'{e} \ggammap-principal  au  sens  de \cite {Ha}.
\end{pr}
     \dem (On utilise les notions et  les  notations  de 
\cite {Ha}.). Soit 
$$\btu \prfib {\wtu }\Gamma '  = \{  (z,(g,w))\; |\;   p'(z)=w,\; 
z\in \btu \}  = \{ (z,(g,p'(z))\; |\; z\in \btu  \}.
$$

     On peut d\'{e}finir une action \`{a} droite de \ggammap  sur 
\btuu  :
$$\btu \prfib {\wtu }\Gamma '\rightarrow \btu
$$ 
$$(z,\gamma ')\rightarrow z\cdot \gamma \stackrel  {def}{=}z\cdot 
g  $$ o\`{u} $\gamma '=(p'(z),g)\in \Gamma '\; $.
( Ici $z\cdot g\; $ est donn\'{e} par l'action naturelle de  
$O(q)\; $  sur le fibr\'{e} $O(q)\; $- principal. )

     On a \'{e}videmment (au cas o\'{u} les produits sont possibles) :
$$     (z\cdot \gamma ')\cdot \gamma '' = z\cdot  (\gamma  '\cdot 
\gamma '') ,$$
$$     z\cdot u = z ,\;  ( u \in \wtu  ),$$
$$     p'(z\cdot \gamma ') = p'(z\cdot g ) = p'(z)\cdot g = 
\alpha '(\gamma ')\; ,\; (\gamma '= (p'(z),g) ).$$

     On  peut  ais\'{e}ment  observer   que \ggammap   op\`{e}re   simplement 
transitivement dans les fibres de $p\; $ .

     Soit ${\cal A} = {(U_i ,\phi _i )_{i\in I}}\; $un atlas  sur 
$M\; $, adapt\'{e}  au feuilletage \mf  ,  tel  que  les  
s\'{e}ctions  locales $s_i\! :\! U\rightarrow \btu $  sont 
feuillet\'{e}s. Ces sections donnent les trivialisations locales  sur 
\btuu : $\; \; h_i'\! :\! p^{-1}(U_i)\rightarrow U_i\times O(q)\; $ 
(v. \cite {Mo}).

     Soit $f_i  = p'\circ s_i\; $  et $U_i\prfib {\wtu }\Gamma  ' 
= \{ (u,(p'(s_i (u)),g)\; |\;  u\in U_i , g\in O(q) \}$.

     En utilisant le diff\'{e}omorphisme naturel
$$U_i\times O(q) \rightarrow  U_i\prfib {\wtu }\Gamma '$$
 on a un diff\'{e}omorphisme
$$h_i\! :\! p^{-1}(U_i )\rightarrow  U_i\prfib {\wtu }\Gamma '$$
qui est $\Gamma '$-\'{e}quivariant et qui se projette sur 
l'identit\'{e} de $U_i\; $ .\qed

\section{La structure diff\'{e}rentielle du groupoide \ggamma}
     Nous allons utiliser ici le groupoide \gabm construit dans 
l'Introduction. Ce groupoide est obtenu comme il suit :

     Soit 
$${\cal H}(\fbtu ) = \{ (z_1 ,z_2 )\; |\;  z_1  ,z_2\in  \btu  \; 
,\;  p'(z_1 ) = p'(z_2 ) \}$$ 
(qui est le produit fibr\'{e} $\btu \prfib {\wtu }\btu \; $  ).
On a une action \`{a} droite de $O(q)\; $ sur ${\cal H}(\fbtu )\; $:
$$ {\cal H}(\fbtu )\times O(q) \rightarrow {\cal H}(\fbtu ))
$$  
$$   ((z_1 ,z_2 ),g) \rightarrow (z_1 ,z_2 )\cdot g\stackrel {def}{=} 
  (z_1 \cdot g , z_2\cdot g)\; $$.

     \ggamma  est  l'espace  des  orbites  de  cette  action.  On  note 
$[z_1 ,z_2]\; $ l'orbite de $(z_1,z_2)\in {\cal H}(\fbtu )\; $  (qui  est  
d\'{e}finie  si  et seulement si $p'(z_1 ) = p'(z_2 ) )$.

     L'inclusion  $I\! :\! M\rightarrow \Gamma \; $ est donn\'{e}e  par  
les  classes diagonales : si $p(z)=x$, alors $I(x) = [z,z]$.

     L'inversion est   donn\'{e}e    par   $[z_1  ,z_2  ]^{-1}  =   
[z_2 ,z_1]\; $  et  les projections \aalpha  et \bbeta  sont  
donn\'{e}es  par  
$\alpha ([z_1 ,z_2 ]) = p(z_2 )\; $  et $\beta ([z_1  ,z_2  ])  = 
p(z_1 )$.

     Nous avons le r\'{e}sultat suivant :
\begin{pr}\label {prd}
Le groupoide \gabm est un groupoide de  Lie 
s\'{e}par\'{e} et compact. Les projections \aalpha et \bbeta  sont  
des  fibrations locales triviales avec la fibre type \ftbuu . 
\end{pr}
     \dem
     Soit ${\cal A}\; $  un  atlas  feuill\'{e}t\'{e}  (consid\'{e}r\'{e}   
d\'{e}j\`{a}   dans   la  d\'{e}monstration  de  {\em la  Proposition 
\ref {pru}} ).  Si  $x\in U_i\; $  ,  alors  soit 
$\overline {F}_{T,i}^1(x)\; $  la fibre de la fibration basique qui  
contienne  $s_i(x)\; $  . 
Pour chaque $i\in I\; $ , la r\'{e}union
$$\Gamma _i= \bigcup _{x\in U_i}  \{s_i (x)\}\times \overline 
{F}_{T,i}^1(x)\subseteq {\cal H}(\fbtu )
$$
est une collection compl\`{e}te de repr\'{e}sentants   pour  les  classes 
d'\'{e}quivalence  qui  constituent $\beta ^{-1}(U_i)\; $.(  La   
projection $\pi '\; $ restriction\'{e}e \`{a} $\Gamma _i\; $  est injective .)

     Soit, pour chaque $i$, une trivialisation locale
$$(p')^{-1}(W_i )\rightarrow W_i\times \ftbu $$  
adapt\'{e}e pour  la  fibration  basique, tel que 
$$f_i (U_i )\subseteq W_i\; ,\; f_i = p'\circ s_i\;  .$$ 
( Pour cela on diminue $U_i\; $ , si c'est  n\'{e}cessaire.  )  On  obtient 
un diff\'{e}omorphisme 
\begin{equation}\label{eqstaru}
\Gamma _i\rightarrow U_i\times \ftbu
\end{equation}
et une bijection
\begin{equation}\label{eqstard}
\beta ^{-1}(U_i)\rightarrow U_i\times \ftbu
\end{equation}

     Si $(s_i (x),y)\in \Gamma _i\; $  et $(s_i(x),z)\in \Gamma _i\; $    
repr\'{e}sentent  la  m\^{e}me  classe,   alors   $y   =   z\cdot 
g_{ij}(x)\; $ , o\`{u}
$$     s_i (x) = g_{ij}(x)\cdot s_j (x)\; ,\; g_{ij}(x)\in O(q)\; .$$

     Les projections \bbeta et \aalpha ont la forme locale
$$\beta (\hat f_x ) = x , \alpha (\hat f_x) = p(f_x ) \; ,$$
o\'{u}  $f_x\in  \overline  {F}^1_{T,i}(x)\;   $et   $\hat f_x\in 
U_i\times \ftbu \; $est obtenu par le  diff\'{e}omorphisme  \rref 
{eqstaru}.

     Il en r\'{e}sulte que $\beta \! :\! \Gamma '\rightarrow  M\; 
$ est  une fibration  locale triviale (avec la fibre type \ftbuu  ) et 
\ggamma est  une   vari\'{e}t\'{e}  s\'{e}par\'{e}e 
avec une trivialisation locale donn\'{e}e par \rref {eqstard}.

     Pour le produit du groupoide \ggamma   on  peut  donner  la  forme 
locale, d'o\'{u} r\'{e}sulte sa diff\'{e}rentiabilit\'{e} :

     Soit $f_x\in \overline {F}^1_{T,i}(x)\; ,\; f_y\in 
\overline {F}^1_{T,i}(y) , p(f_x ) = y  =  p(s_j (y))$.  On 
r\'{e}sulte    qu'il    existe    $g\in    O(q)\;    $tel     que 
$f_x=s_j(y)\cdot g\; $  .  Soit 
$f'_x = f_y\cdot g\; $; alors  $\hat  f_x\cdot  \hat  f_y   =\hat 
f'_x\; $ .

     En utilisant la trivialisation \rref {eqstaru} pour la forme
 locale de $I\; $ on peut constater que $I(M)\subset \Gamma \; $ est une 
sous-vari\'{e}t\'{e}. Nous  allons identifier les points de  $M\; 
$ avec leurs images par $I\; $.

     Pour trouver la forme locale de  l'inversion  du   groupoide 
\ggamma  soit $f_x\in \overline {F}^1_{T,i}(x)\; ,\; y=p(f_x)  =  
p(s_j (y))$.  Il  r\'{e}sulte  qu'il 
existe $g\in O(q)\; $ tel que $f_x =s_j(y)\cdot g\; $. Soit  $f_y  
= s_i(x)\cdot g^{-1}\in \overline {F}^1_{T,j}(y)$.
Alors  la  forme  locale   de   l'inverse    fait    correspondre  
$\hat f_y\;  $ \`{a}  $\hat f_x\; $  d'o\'{u} il r\'{e}sulte que 
l'inversion est diff\'{e}rentiable. Elle est aussi  un 
diff\'{e}omorphisme parce qu'elle est bijective et involutive.

     L'inversion   transforme    $\alpha $-fibres   en $\beta $-fibres    (et 
r\'{e}ciproquement)  et      laisse       invariante       la  
sous-vari\'{e}t\'{e}  $M\; $.  Il 
r\'{e}sulte que $\alpha \! :\! \Gamma \rightarrow M\; $ est aussi une  
fibration   avec   la m\^{e}me fibre type \ftbuu  . \qed

     Si on utilise les constructions faites dans la d\'{e}monstration 
de {\em la Proposition \ref {prd}}, on peut voir qu'il existe une 
trivialisation locale pour \hfbtuu de telle facon qu'on a le 
r\'{e}sultat suivant :
\begin{pr}\label {prt} $\pi \! :\! \hfbtu \rightarrow  \Gamma  \; 
$est un fibr\'{e} $O(q)$-principal.
\end{pr}

     L'orbite du groupoide \gabm qui passe par $u\in M\; $ est
$$\alpha (\beta ^{-1}u)) = \beta (\alpha ^{-1} (u))$$
et  l'ensemble  des  orbites  constitue  un    feuilletage    (en  
g\'{e}n\'{e}ral singulier, de Stefan) (cf. \cite {AD} ).

\begin {pr}\label {prp} Les orbites du groupoide \ggamma sont  les feuilles 
du feuilletage \mfb .
\end{pr}
\dem     Soit $u\in M\; $ et $\overline {F}\; $ la feuille du feuilletage  
\mfb, $u\in \overline {F}$. 
Si $z\in p^{-1}(u)\; $ et \ftbuu  est la feuille du feuilletage 
$(\btu ,\ftbu )\; $qui contient $z\; $ , alors  la  restriction  
$p\! :\! \ftbu \rightarrow \overline {F}\; $  
est  une fibration principale \cite {Mo} et en tenant compte de 
la d\'{e}monstration de {\em la  Proposition  \ref  {prd}}  ,  il 
r\'{e}sulte (pour un cerain $z$) :
$$\alpha (\beta ^{-1}(u))=p(\ftbu )= \overline {F}
$$
d'o\'{u} le r\'{e}sultat. \qed 

\section {L'\'{e}quivalence des groupoides \ggamma et \ggammap}
     Le  but  de   ce   paragraph   est   de   d\'{e}montrer   le  
r\'{e}sultat 
suivant:
\begin{pr} \label {prc}
\btuu  est une \'{e}quivalence de \ggamma sur \ggammap  (au  sens 
de \cite {Ha}).\end{pr}
     La d\'{e}monstration de cette proposition sera donn\'{e}e  par  deux 
lemmes.
\begin{lm}\label {lmunu}
\btuu  est un morphisme (g\'{e}n\'{e}ralis\'{e}) de \ggamma  dans 
\ggammap  (au sens de \cite {Ha}).\end{lm}
\dem Soit $p\! :\! \btu \rightarrow M\;  $   et  $p'\!  :\!  \btu 
\rightarrow \wtu $.  Nous  avons d\'{e}montr\'{e} 
({\em Proposition 1}) le fait que \btuu est un fibr\'{e} 
$\Gamma '$-principal et  donc \ggammap op\`{e}re \`{a} droite sur  
\btuu relativement \`{a} $p'$.

     Pour d\'{e}montrer la {\em Lemme}, on doit d\'{e}finir une  action 
\`{a} gauche de \ggamma sur \btuu , relativement  \`{a}  $p$,  tel  que  
les  deux  actions commutent (cf. \cite {Ha}).

     Soit $\gamma = [z_1 ,z_2 ]\in \Gamma ,\; \alpha (\gamma ) = p(z_2) = x , 
\; \beta (\gamma ) = p(z_1) = y\; $et  $z\in  \btu  ,\;   p(z)  = 
\beta(\gamma )=y\; $, donc il  existe  $g\in O(q)\; $  tel  que  
$z_1=  z\cdot g\; $. On d\'{e}finit $\gamma \cdot z = z_2\cdot g\; $.
     On peut v\'{e}rifier directement les \'{e}galit\'{e}s :
$$\gamma _1\cdot (\gamma _2\cdot z)=(\gamma _2\cdot \gamma _2)\cdot z,\; 
\gamma _1,\gamma _2\in \Gamma ,z\in \btu ,\; \alpha (\gamma _1 )  = 
\beta (\gamma _2) ,\; p(z) = \alpha (\gamma _2 ) ;$$
$$ p(\gamma \cdot z) = \beta (\gamma ) ,\; \gamma \in \Gamma  ,\; 
z\in \btu,\; ,\; p(z) = \alpha (\gamma ) ;$$
$$u\cdot z = z ,\; u\in M\subset \Gamma , z\in \btu ,\;  p(z) = u 
$$et le fait que les deux actions commutent. \qed

     Soit (cf. \cite {Ha}) $(\btu )^0\; $ la m\^{e}me vari\'{e}t\'{e}   
\btuu munie  des  m\^{e}mes applications $p\; $ et $p'\; $ dans 
$M\; $  et \wtuu  mais cette fois \ggamma et \ggammap op\'{e}rant 
\`{a} droite et \`{a} gauche respectivement en posant
$$z\circ \gamma = \gamma ^{-1}\cdot z\; \; ,\; \; 
    \gamma '\circ z = z \cdot (\gamma ')^{-1} .$$
\begin {lm} $(\btu )^0\; $  est un fibr\'{e} $\Gamma  $-principal 
sur \wtuu de projection $p'\! :\! \btu \rightarrow \wtu \; $.\end{lm}
\dem  On a le  produit  fibr\'{e}  des  applications  $p\;  $  et 
$\gamma\; $ :
$$\btu \prfib {M}\Gamma = \{ (z,\gamma )\in \btu \times \Gamma \; 
|\;  p(z) = \beta (\gamma ) \}$$
     L'action \`{a} droite de \ggamma sur $(\btu )^0\; $  est :

     -Simple : Soit $z\circ \gamma _1= z\circ \gamma _2\; 
(\gamma _1 = [z_1 ,z'_1], \gamma _2= [z_2 ,z'_2] )\; $
ou $[z'_1,z_1 ]\cdot z = [z'_2,z_2 ]\cdot z\; $.  On  a  $p(z)  = 
p(z'_1) = p(z'_2),\;  z'_1\cdot g_1= z = z'_2\cdot g_2\; $  et 
$z_1\cdot g_1 = z_2\cdot g_2\; (g_1 , g_2\in  O(q)  )\;  $   donc 
$\gamma _1=\gamma _2$  .

     -Transitive sur les fibres de $p'\; $: Soit $z, z'\in \btu 
\; $   tel  que $p'(z) = p'(z')$. Si $g = [z',z]$, alors  $z\circ 
\gamma = z'\; $ .

     Soit $t_{\alpha }\! :\! W_{\alpha }\rightarrow \btu \; $
les   sections  locales   de   la fibr\'{e} 
locale triviale $p'\! :\! \btu \rightarrow \wtu \; $  (  avec  la  
fibre  type \ftbuu ) donn\'{e}es par les trivialisations locales
$$(p')^{-1}  (W_{\alpha } )\simeq \wtu \times \ftbu $$

     On a
$$W_{\alpha }\prfib {M}\Gamma =  \{  (w,\gamma )\in  W_{\alpha  }\times 
\Gamma \; |\;  p(t_{\alpha }(w))  = \beta (\gamma ) \} = 
     = (p\circ t_{\alpha } )^{\ast }(\Gamma )$$
o\'{u}  la  derni\`{e}re   expression   repr\'{e}sente   l'image   
r\'{e}ciproque locale de la fibration  locale   triviale  $\Gamma 
\stackrel {\beta }{\rightarrow }M\; $qui a la m\^{e}me fibre  type 
\ftbuu .  Pour $W_{\alpha }\; $   suffisamment  diminu\'{e},   on   
a   un diff\'{e}omorphisme
$$h_{\alpha  }  \!  :\!  (p')^{-1}(W_{\alpha  }  )    \rightarrow  
W_{\alpha }\prfib {M}\Gamma $$
qui se projette sur l'identit\'{e} de $W_{\alpha }$. Par un  
calcul  direct  on peut constater aussi la $\Gamma $-\'{e}quivariance de  
$h{_\alpha }\; $( $h{_\alpha } (z\circ \gamma  )=  h_{\alpha } (z)\circ 
\gamma \; $, o\'{u} la derni\`{e}re action est donn\'{e}e par 
$(w,\gamma _1)\circ \gamma =(w,\gamma _1\circ \gamma )\; $ ).\qed 

\section {La structure m\'{e}trique.}
     Nous allons utiliser ici une m\'{e}trique  riemannienne  sur 
\btuu  donn\'{e}e par son parall\'{e}lisme transvers,  determin\'{e}  
par une connexion m\'{e}trique adapt\'{e}e au feilletage \mf ( cf 
\cite  {Mo}).  Les   r\'{e}sultats   cuntenus   dans   les   {\em 
Lemme \ref {ltrei}, Propositions \ref {prss}, \ref {prsp} et 
\ref {pro} } sont 
connus, mais il est important de pr\'{e}ciser  les  details  pour 
construir dans  {\em  la  Proposition  \ref {prn}} la  structure 
m\'{e}trique du \ggamma .

\begin{pr}\label {prss}
 Soit \mf un feuilletage riemannien et $g\; $ une m\'{e}trique 
riemannienne sur  $M\; $  adapt\'{e}e  au feuilletage.  Alors  il 
existe une m\'{e}trique riemannienne $g\; $ sur \btuu  tel que :

     a)  La  projection $p\! :\! \btu \rightarrow M\; $est  une   
submersion riemannienne.

     b) Pour $a\in O(q)\; $ la translation $R_a\; $est une isom\'{e}trie.

     c)  La   m\'{e}trique $\stackrel {\sim }{g}\; $est  une  m\'{e}trique  
adapt\'{e}   pour   le feuilletage \'{e}lev\'{e} $(\btu ,\frtu )\; $ qui est 
riemannien.
\end{pr}     
\dem Soit $\omega \! :\! T\btu \rightarrow o(q)\; $la  forme  de  connexion  
de  la connexion de Levi Civita associ\'{e} \`{a} la m\'{e}trique 
$\stackrel {\sim }{g}\; $ et  $<\; ,\; >\; $  la m\'{e}trique de 
Killing associ\'{e} \`{a} l'alg\`{e}bre le Lie $o(q)\; $. On 
d\'{e}finit :
$$\stackrel {\sim }{g}(\stackrel {\sim }{X},\stackrel {\sim }{Y}) = 
g(p_{*}\stackrel {\sim }{X} , p_{*}\stackrel {\sim }{Y}) + 
< \omega (\stackrel {\sim }{X}), \omega  (\stackrel  {\sim  }{Y}) 
>$$

     L'affirmation {\em a)} d\'{e}coule directement de la d\'{e}finition .

     Pour {\em b)}  on  utilise  l'invariance  de  la  m\'{e}trique  de 
Killing  \`{a} l'adjunction.

     Pour {\em c)} on peut partir du fait que la forme $\omega \; $ 
est localement projectable.   La   distribution   verticale   est    
orthogonale  \`{a} la distribution donn\'{e}e par les feuilles de 
\frtuu  et le produit de  deux champs  (locals) basiques est  une 
fonction  basique ( pour \frtuu ).\qed 

     Cette m\'{e}trique sera  utilis\'{e}e   pour   induire   des  
m\'{e}triques riemanniennes sur $\wtu ,\; {\cal H}(\fbtu )\;  $et 
\ggamma .
\begin{pr}\label {prsp} Si $\stackrel {\sim}{g}\; $ est une 
m\'{e}trique adapt\'{e} au  feuilletage riemannien 
\bff , alors il existe une m\'{e}trique riemannienne $g'\; $ 
sur \wtuu , tel que la submersion $p'\! :\! \btu \rightarrow \wtu \; $ 
 soit riemannienne.
\end{pr}     
\dem Le feuilletage \bff a l'adh\'{e}rence r\'{e}guli\`{e}re et 
on  peut appliquer  {\em Lemme 5.2  (\cite {Mo},  pg.  156)}.  Il  
r\'{e}sulte   que   le feuilletage \bffb est riemannien  ,  il  
a  comme  m\'{e}trique  adapt\'{e}e  la  m\'{e}trique  $\stackrel 
{\sim }{g}\; $  et sa structure transverse d\'{e}pend  seulement de la 
structure transverse du feuilletage \bff . On  peut  projeter 
donc la m\'{e}trique $\stackrel {\sim}{g}\; $ sur une m\'{e}trique 
$g''\; $ sur \wtuu  de telle fa\c{c}on que la submersion $p'\; $ soit 
riemannienne. \qed 

     Dans ce qui suit, nous allons consid\'{e}rer la m\'{e}trique  
$\stackrel {\sim }{g}\; $   sur  \btuu  donn\'{e}e   par   {\em   la 
Proposition \ref {prss}} .
\begin{pr}\label {pro}Le goupoide transformationnel 
$$\Gamma ' = \grupoid {\wtu \times O(q)}{\wtu }{\alpha '}{\beta '}$$
est riemannien.
\end{pr} 
\dem On peut d\'{e}finir sur \ggammap la  m\'{e}trique  produit  de  la 
m\'{e}trique $g'\; $sur \wtuu (donn\'{e}e par la {\em Proposition 
\ref {prsp}}) avec la m\'{e}trique  de  Killing sur $O(q)\; $. En  tenant  
compte  des  expressions  qui  donnent  les 
projections $\alpha ' ,\; \beta '\; $ et l'inversion du groupoide 
\ggammap , il r\'{e}sulte  que $\beta '\; $  est  une submersion 
riemannienne; pour d\'{e}montrer que $\alpha '\; $ est aussi  une 
submersion riemannienne et  l'inversion  est une isom\'{e}trie, 
il est suffisant de d\'{e}montrer  que  l'action   de  $O(q)\;  $ 
sur \wtuu est 
faite par des isom\'{e}tries.  La  conclusion  de  la 
{\em Proposition} r\'{e}sulte donc du {\em Lemme} suivante :
\begin{lm}\label {ltrei}
  Soit $g\in O(q),\; z\in \btu ,\; V\; $un champ  vectoriel 
feuill\'{e}t\'{e}  pour \bffb , orthogonal sur  les  feuilles  et 
$\stackrel {\sim }{V}\; $le champ vectoriel proj\'{e}t\'{e} de $V\; $ 
(sur \wtuu ).Alors $R_{g*}(V)\; $ est orthogonal sur les feuilles 
et
$$ ||R_{g*}(v)_{z\cdot g}||_{\btu }  = || V_z ||_{\btu } \; ;$$
$$ ||R_{g*}(\stackrel {\sim }{V})_{p'(z\cdot g)}||_{\wtu }= 
||\stackrel {\sim }{V}_{p'(z)}||_{\wtu }\; .$$
\end{lm}     
\dem  $r_{g*}\; $  est une application  feuill\'{e}t\'{e}   pour   
\bffb  (\cite {Mo}).  Il r\'{e}sulte de la {\em Proposition  \ref 
{prss}} qu'elle est aussi une isom\'{e}trie, donc 
la premi\`{e}re relation est d\'{e}montr\'{e}e.

     De la {\em Proposition \ref {prsp}} il r\'{e}sulte que $p'\; $  est  
une  submersion riemannienne. Les vecteurs $R_{g*}(v)_{z\cdot g}\; 
$ et $V_z\; $sont orthogonaux  sur les fibres de $p'\; $ et ils se 
projettent par $p'\; $  sur  les  vecteurs 
$R_{g*}(\stackrel {\sim }{V})_{p'(z\cdot g)}\; $ et 
$\stackrel {\sim }{V}_{p'(z)}\; $respectivement, d'o\'{u}     le 
r\'{e}sultat.\qed 

     Le r\'{e}sultat suivant va pr\'{e}ciser la structure m\'{e}trique 
de \ggamma :
\begin{pr}\label {prn} $\pi \! :\! \hfbtu \rightarrow \Gamma \; $
  est  un morphisme  de groupoides de Lie. Il existe sur \ggamma une 
m\'{e}trique riemannienne  tel que $\pi \; $ soit une submersion 
riemannienne et le groupoide \gabm soit riemannien.
\end{pr} 
\dem  Le feuilletage \bffb est simple et a  donc  une  holonomie 
nulle (\cite {Mo}). Son groupoide d'holonomie
$$\grupoid {\hfbtu }{\btu }{p_1}{p_2}
$$
est s\'{e}par\'{e} , il est un groupoide riemannien  au  sens  de 
\cite {GGHR}  et sa m\'{e}trique est donn\'{e}e par
$$||V||_{\cal H}^z = ||p_{1*}(V)||_{\btu }^{z_1}+ 
||p_{2*}(V)||_{\btu }^{z_1} - ||(p'\circ p)_*(V)||_{\btu }^{z_1}
$$(cf. \cite {Wi}).

     Les  actions  \`{a}  droite  de  $O(q)\; $  sur  \hfbtuu  et 
\btuu  commutent avec $p_1\; $  et $p_2\; $ . En  tenant  compte   que   les

actions  \`{a} droite de $O(q)\; $ sur \btuu  et \wtuu sont des 
isom\'{e}tries , il r\'{e}sulte  que l'action \`{a} droite de 
$O(q)\; $ sur \hfbtuu  est  aussi  faite  par  des 
isom\'{e}tries. On peut d\'{e}finir donc  une  m\'{e}trique  sur  
\ggamma  tel  que $\pi \; $  soit  une  submersion  riemannienne.

     Les diagrammes suivants sont commutatives :
\newline
\begin{picture}(240,65)(-47,5) 
\put  (6,10){$B_T^1$}  \put  (0,50){${\cal  H}(\overline  {\cal 
F}_T^1)$} \put (64,10){$M$} \put (66,50){$\Gamma $} 
\put (14,46){\vector (0,-1){26}} \put (68,46){\vector (0,-1){26}} 
\put (22,14){\vector (1,0){39}} \put (35,54){\vector (1,0){28}} 
\put (40,18){$p$}  \put  (40,58){$\pi $}  \put  (16,30){$p_1$}  \put 
(70,30){$\alpha $} 
\put  (156,10){$B_T^1$}  \put  (150,50){${\cal  H}(\overline  {\cal 
F}_T^1)$} \put (214,10){$M$} \put (216,50){$\Gamma $} 
\put (164,46){\vector (0,-1){26}} \put (218,46){\vector (0,-1){26}} 
\put (172,14){\vector (1,0){39}} \put (185,54){\vector (1,0){28}} 
\put (190,18){$p$} \put  (190,58){$\pi $}  \put  (166,30){$p_2$}  \put 
(220,30){$\beta $} 
\end{picture}

\begin{picture}(240,65)(-30,5) 
\put  (9,10){$\Gamma $}  \put  (0,50){${\cal  H}(\overline  {\cal 
F}_T^1)$} \put (63,10){$\Gamma $} \put (52,50){${\cal  
H}(\overline  {\cal F}_T^1)$} 
\put (14,46){\vector (0,-1){26}} \put (67,46){\vector (0,-1){26}} 
\put (18,14){\vector (1,0){43}} \put (35,54){\vector (1,0){18}} 
\put   (40,18){$*^{-1}$}    \put    (40,58){$*^{-1}   $}     \put  
(16,30){$\pi $}  \put (69,30){$\pi $} 
\put  (136,10){$\Gamma \times \Gamma $}  
\put  (110,50){${\cal  H}(\overline {\cal F}_T^1)\times 
{\cal  H}(\overline {\cal F}_T^1)$} 
\put (224,10){$\Gamma $} \put (212,50){${\cal  H}(\overline {\cal 
F}_T^1)$} 
\put (150,46){\vector (0,-1){26}} \put (227,46){\vector (0,-1){26}} 
\put (165,14){\vector (1,0){56}} \put (192,54){\vector (1,0){20}} 
\put (195,18){$\mu $} \put  (195,58){$\mu '$}   \put   (152,30){$\pi 
\times \pi $}  \put (230,30){$\pi $} 

\end{picture}
\newline
(o\'{u} on a not\'{e} par $*^{-1}\; $l'inverse et par 
$\mu ,\; \mu '\; $les produits  des groupoides ). Il r\'{e}sulte le  
fait  que  $\pi \; $  est  un  morphisme des groupoides de Lie.

     Le fait que l'inversion du groupoide \ggamma  est  une  
isom\'{e}trie r\'{e}sulte de  la  construction  de  m\'{e}trique  
de  \ggamma .  Pour  que  le groupoide \ggamma soit riemannien il reste 
\`{a}  d\'{e}montrer que \aalpha  est  une submersion riemannienne .

     Dans le premier diagramme ci-dessus,  les $\alpha $-fibres sont  les 
projections des $\pi $-fibres. Soit $\gamma \in \Gamma \; $ et 
$V\in T_{\gamma}(\gamma  )\;  $  orthogonal  \`{a}  la  fibre  de 
\aalpha qui contient $\gamma \; $.  Soit  $h\in  \hfbtu  ,\;  \pi 
(h)=\gamma \; $et $W\in T_h\hfbtu \; $orthogonal  \`{a}  la  $\pi 
$-fibre, qui se projette sur $V\; $ et qui a  la  norme  \'{e}gale  
\`{a}  celle   de   $V\; $.   En   utilisant   la  
commutativit\'{e}  de la  diagramme  consid\'{e}r\'{e}e,  il  
r\'{e}sulte  que $W\; $ est  orthogonal  \`{a}  la 
$p_1$ - fibre et se projette par $p_{1*}\; $sur un  vecteur $W'\; 
$  de $T_z\btu \;,\; z=p_1(h)\; $ , qui a  la  m\^{e}me  norme  
que $W\; $.  Mais $W\; $ est  aussi orthogonal \`{a} 
la fibre de $p\; $ et il se projette par $p_{1*}$  sur  un  vecteur  
$W'\in T_z\btu\; ,\; z=p_1(h)\; $,  qui  a  la   m\^{e}me   norme  
que $W\; $.  Mais $W'\; $  est  aussi orthogonal \`{a} la fibre 
de $p\; $ et il  se  projette  par $p_*\; $    sur  le vecteur  
$p_*\circ p_{1*}(W) = \alpha _*\circ \pi _* (W) = \alpha _* (V)\; 
$ qui a la m\^{e}me norme  que  $V\; $. Donc \aalpha est une submersion 
riemannienne.\qed 
\vspace {1 cm}

     Note. {\em L'auteur remercie M. le Professeur Gilbert  Hector  de 
la  patience  manifest\'{e}  durant  l'initiation   \`{a}   l'\'{e}tude   des 
Groupoides Riemanniens et d'avoir indiqu\'{e} ce probl\`{e}me et 
M. le Professeur Mircea Puta pour les conversations tr\'es utiles sur 
ce sujet.}


\end{document}